

 \input amstex
 \documentstyle{amsppt}
 \magnification=\magstep1

\TagsOnRight

\hyphenation{auto-mor-phism auto-mor-phisms co-homo-log-i-cal co-homo-logy
co-homo-logous dual-izing pre-dual-izing geo-metric geo-metries geo-metry
half-space homeo-mor-phic homeo-mor-phism homeo-mor-phisms homo-log-i-cal
homo-logy homo-logous homo-mor-phism homo-mor-phisms hyper-plane hyper-planes
hyper-sur-face hyper-sur-faces idem-potent iso-mor-phism iso-mor-phisms
multi-plic-a-tion nil-potent poly-nomial priori rami-fication sin-gu-lar-ities
sub-vari-eties sub-vari-ety trans-form-a-tion trans-form-a-tions Castel-nuovo
Enri-ques Lo-ba-chev-sky Theo-rem}

\define\red{_{\roman{red}}}              
\define\mi{_{\roman{min}}}               
\define\rest#1{_{\textstyle{\vert}#1}}   
\define\iso{\cong}
\define\Span#1{\left<#1\right>}          
\define\na{\nabla}
\define\nao{\na_0}
\define\half{{\textstyle{1\over2}}}
\define\Vol{\mathop{\roman{vol}}}        
\define\dd{\roman{d}}                    

\define\C{\Bbb C}    
\define\R{\Bbb R}    
\define\Z{\Bbb Z}    
\define\proj{\Bbb P} 
\define\plane{\Bbb{CP}^2} 

\define\sA{{\Cal A}} 
\define\sB{{\Cal B}} 
\define\sG{{\Cal G}} 
\define\sJ{{\Cal J}} 
\define\sM{{\Cal M}} 
\define\Oh{{\Cal O}} 
\define\sS{{\Cal S}}

\define\al{\alpha}
\define\de{\delta}
\define\ep{\varepsilon}
\define\fie{\varphi}
\define\ga{\gamma}
\define\ka{\kappa}
\define\om{\omega}
\define\si{\sigma}
\define\De{\Delta}
\define\Ga{\Gamma}

\define\Om{\Omega}

\define\ad{\operatorname{ad}}          
\define\alg{\operatorname{alg}}
\define\card{\mathop{\#}}              
\define\codim{\mathop{{\roman{codim}}}\nolimits} 
\define\coker{\mathop{{\roman{coker}}}} 
\define\diff{\operatorname{diff}}
\define\rank{\mathop{{\roman{rank}}}}  
\define\var{\operatorname{var}}
\define\vcodim{\operatorname{v.codim}}
\define\vdim{\operatorname{v.dim}}
\redefine\mod{\mathop{{\roman{mod}}}}  
\define\Bs{\mathop{{\roman{Bs}}}}      
\define\Diff{\operatorname{Diff}}
\define\Grass{\operatorname{Grass}}    
\define\Hilb{\mathop{{\roman{Hilb}}}\nolimits} 
\define\Hom{\mathop{{\roman{Hom}}}\nolimits}   
\define\Mod{\operatorname{Mod}}
\define\Mon{\operatorname{Mon}}
\define\Pic{\operatorname{Pic}}        
\define\Sing{\mathop{{\roman{Sing}}}}  
\define\W{\mathop{\bigwedge^2}}        

\redefine\O{\operatorname{O}}
\define\SO{\operatorname{SO}}
\define\SU{\operatorname{SU}}
\define\Spin{\operatorname{Spin}}
\define\U{\operatorname U}             

 \document

 \topmatter
 \title Spin canonical invariants of\\ 4-manifolds and algebraic
surfaces\endtitle
 \author Andrei Tyurin \endauthor
 \affil Steklov Math Institute \endaffil
 \address Algebra Section, Steklov Math Institute, Ul\. Vavilova 42, Moscow
117966 GSP--1, Russia \endaddress
 \email Tyurin\@top.mian.su {\it or} Tyurin\@CFGauss.Uni-Math.Gwdg.De
{\it or}\newline Tyurin\@Maths.Warwick.Ac.UK
\endemail

 \endtopmatter
 \rightheadtext{Spin canonical invariants of 4-manifolds}

\head 0. Introduction\endhead

Donaldson's celebrated construction of polynomial invariants for 4-manifolds
(see for example \cite{D--K}) has served as a model for an enormous range of
recent results and constructions relating algebraic geometry and differential
topology. The short written version of Donaldson's talk at the 1990
Arbeitstagung contained the construction of jumping instantons, and was the
starting point for our constructions and applications of the new spin
polynomial
invariants of 4-manifolds and algebraic surfaces. I am very grateful to the
Royal Society of London for the opportunity to make an extremely pleasant and
stimulating 3 months visit to Warwick, Oxford and Cambridge in summer 1993,
which has allowed me to discuss these constructions with Simon Donaldson.

I have also benefited from several conversations with Nigel Hitchin. His
questions about what Yang--Mills connections look like from a geometric point
of
view (see for example \cite{H}) has involved a journey to classical incidence
geometry, or ``algebraic protogeometry''.

To explain what I mean by this, suppose that we want to prove the following
differentiable version of the Poincar\'e conjecture for $\plane$:

\proclaim{Conjecture \rm(DPC for $\plane$)} The complex projective plane
$\plane$ has a unique differentiable structure. \endproclaim

To prove this, we have to say what is a {\it line} or a {\it nonsingular conic}
of $\plane$ in terms of the underlying differentiable structure of $\plane$. If
we can do this, and check that the lines satisfy the classical incidence
axioms,
then we are done. If this can't be done, then equipping $\plane$ with any
Riemannian metric $g$, we have to:

\roster
\item say what it means in terms of Riemannian geometry; and

\item describe the dependence on the metric $g$.
\endroster

Now the analog of a {\it nonsingular conic} in the language of Riemannian
geometry is a $g$@-instanton of topological type $(2,0,2)$ up to gauge
equivalence, that is, a $g$@-antiselfdual $\SU(2)$@-connections on a vector
bundle $E$ with $c_2=2$. This holds because of Donaldson's identification
 $$
\left\{
\matrix
\text{instantons for the Fubini--Study}\\
\text{metric on $\plane$}
\endmatrix
\right\}
=
\bigl\{\text{stable holomorphic vector bundles}\bigr\},
 $$
and Barth's interpretation of stable bundles on $\plane$ (see \cite{B}).

The next question is:
 $$
\text{How many nonsingular conics can be inscribed in a general 5-gon?}
 $$
This constant is the Donaldson polynomial $\ga^5_{\plane}$ evaluated at the
generator of $H^2(\plane,\Z)$. Thus the Donaldson polynomials provide the
information about the incidence correspondence between conics and lines.

My task today is to introduce new invariants of this geometric type. I avoid
important technical details, referring the reader to more specialised articles,
but my aim is to indicate {\it why the following facts are actually true}.

 \head Contents \endhead

1. Construction of spin polynomials invariants

2. Dependence on the metric

3. How much is the independence?

4. Applications

 \head 1. Construction of spin polynomials invariants \endhead

For a smooth, compact 4-manifold $M$, the Stiefel--Whitney class is the
characteristic vector $w_2(M)\in H^2(M,\Z/2)$ of the intersection form $q_M$
on $H^2(M,\Z)$. Any vector $C\in H^2(M,\Z)$ with $C\equiv w_2(M)\mod2$ is
called a {\it $\Spin^\C$ structure} of $M$. Thus a $\Spin^\C$ structure on
$M$ is just a lifting of $w_2$ to an integer class.

For the rest of this paper, for $\si\in H^2(M,\Z)$, we write $L_\si$ to
denote a complex line bundle with first Chern class
 $$
c_1(L_\si)=\si.
 $$
Any Riemannian metric $g$ and $\Spin^\C$ structure $C$ on $M$ defines a
decomposition of the complexified tangent bundle $TM_\C$ as a tensor product
 $$
TM_\C=(W ^-)^*\otimes W^+
 $$
of two rank 2 Hermitian vector bundles $W^\pm$ with
 $$
\W W^\pm=L_C.
 $$

Moreover, for any $\U(2)$@-bundle $E$ on $M$ and any Hermitian connection
$a\in{\sA}_h$ on $E$, putting any Hermitian connection $\nao$ on $L_C$ gives a
coupled Dirac operator
 $$
D^{g,C,\nao}_a\:\Ga^\infty(E\otimes W^+)\to\Ga^\infty(E\otimes W^-).
 \tag1.1
 $$
Now the orbit space of irreducible connections modulo the gauge group
 $$
\sB(E)=\sA^*_h(E)/\sG
 $$
contains the subspace
 $$
\sM^g(E)\subset\sB(E)
 \tag1.2
 $$
of antiselfdual connections with respect to the Riemannian metric $g$.

For any positive integer $r$, we can consider the subspace of
jumping bundles:
 $$
\sM^{g,C}_r(E)=
\bigl\{(a)\in\sM^g(E)\,\,\bigm|\,\,\rank\ker D^{g,C,\nao}_a\ge r\bigr\}
\subseteq\sM^g(E).
 \tag1.3
 $$
The number $r$ is called the {\it jumping level} of $E$. The collection of
these
subspaces defines a filtration:
 $$
\sM^g(E)\supseteq\sM_1^{g,C}(E)\supseteq\cdots
\supseteq\sM_r^{g,C}(E)\supseteq\cdots.
 \tag1.4
 $$
The virtual (expected) codimension of $\sM^{g,C}_r(E)$ is given by
 $$
\vcodim\sM^{g,C}_r(E)=2r^2-2r\chi_C(E),
 \tag1.5
 $$
where $\chi_C(E)$ is the index of the coupled Dirac operator (1.1), which
depends only on the Chern classes of $E$ and the $\Spin^\C$ structure $C$.


The analog of the Freed--Uhlenbeck theorem, that for generic metric $g$ the
moduli space $\sM^g(E) $ (2.1) is a smooth manifold of the expected dimension
with regular ends (see \cite{F--U}, Theorem~3.13), was proved in \cite{P--T},
Chap\.~2, \S3 for the first step $\sM_1^{g,C}(E)$ of the filtration (1.4).
Moreover, $\sM^g(E)$ admits a natural orientation (see \cite{D--K}) inducing an
orientation on $\sM^{g,C} _r(E)$, because its normal bundle has a natural
complex structure. This orientation is described in details in \cite{P--T},
Chap\.~1, \S5.

We need the usual restrictions on the topology of $M$: we suppose that
 $$
b^+_2(M)=2p_g(M)+1
$$
is odd. Then both
 $$
\vdim\sM^g(E)=2d\qquad\text{and}\qquad\vdim\sM^{g,C}_r(E)=2d_r
 $$
must be even.

It is very natural here to use the Uhlenbeck compactification of the first
step of our filtration. We get a filtration
 $$
\overline{\sM^g(E)}\supseteq\overline{\sM_1^{g,C}(E)}
\supseteq\cdots\supseteq\overline{\sM_r^{g,C}(E)}\supseteq\cdots
 \tag1.4$'$
 $$
Actually, in applications, when $r$ is big (say $\ge5$), we only consider cases
with compact moduli spaces $\sM_r^{g,C}(E)$.

Now for any element of our filtration and for a general metric $g$, slant
product defines a cohomological correspondence
 $$
\mu\:H_i(M,\Z)\to H^{4-i}(\overline{\sM_r^{g,C}(E)},\Z),
 $$
and a collection of polynomials
 $$
\ga^{E}_g,s\ga^{E}_{g,C},\dots,s_r\ga^{E}_{g,C,},\ldots\in S^*H^2(M,\Z),
 \tag1.6
 $$
where the first $\ga^{E}_g$ is the Donaldson polynomial, the second
$s\ga^{E}_{g,C}$ is called the {\it spin polynomial}, and the general term the
{\it spin polynomial of jumping level $r$}.

Since
 $$
\overline{\sM^g(E \otimes L)}=\overline{\sM^g(E)}
\qquad\text{for any line bundle $L$,}
 $$
it follows that $\overline{\sM^g(E)}$ depends only on $c_1\mod2$ and on the
first Pontryagin number $p_1=c_1^2-4c_2$, where $c_i$ are the Chern classes of
$E$. Hence the Donaldson polynomial $\ga^{E}_g$ depends only on $c_1\mod 2$ and
$p_1$:
 $$
\ga^{E}_g=\ga^{c_1\mod2,p_1}_g.
 $$

Also, for the twisted bundle $E\otimes L$, the coupled Dirac operator
 $$
D^{g,C-2c_1(L),\nao}_{a\otimes L}\:\Ga^\infty(E\otimes L\otimes W^+\otimes L^*)
\to\Ga^\infty(E\otimes L\otimes W^-\otimes L^*)
 $$
for the shifted $\Spin^\C$ structure $C-2c_1(L)$ is precisely the Dirac
operator (1.1). Hence the space $\overline{\sM_r^{g,C}(E)}$ depends only on
$c_1+C\in H^2(M,\Z)$ and $p_1$. Also, the spin polynomial $s_r\ga^{E}_{g,C,}$
of
jumping level $r$ depends only on $c_1+C$ and $p_1$:
 $$
s_r\ga^{E}_{g,C}=s_r\ga^{p_1}_{g,C+c_1}.
 \tag1.7
 $$

\example{Example: Algebraic surfaces} If $M$ is the underlying manifold of an
algebraic surface $S$, then there exists a canonical $\Spin^\C$ structure given
by the anticanonical class $-K_S$ (we drop the lower index when there is no
danger of confusion).

 In this case, for any $H\in\Pic S\subset H^2(S,\Z)$, by the
Donaldson--Uhlenbeck identification theorem, we have
 $$
\sM^{g_H}(E)=M^H(2,c_1,c_2),
 $$
where the right-hand side is the moduli space of holomorphic $H$@-slope stable
bundles on $S$ with Chern classes $c_1,c_2$.

Making this identification $(a)=E$, we have identifications
 $$
\ker D^{g_H,-K}_a=H^0(E)\oplus H^2(E)
\qquad\text{and}\qquad
\coker D^{g_H,-K}_a=H^1(E),
 \tag1.8
$$
where $H^i(E)$ denote coherent cohomology groups. Thus the index of the coupled
Dirac operator for this special $\Spin^\C$ structure is
 $$
\chi_{-K_S} D^{g_H,-K}_a=h^0(E)-h^1(E)+h^2(E).
 $$
So the subspace $\sM_r^{g_H,-K_S}(E)$ is the {\it Brill--Noether} locus
 $$
\sM^{g_H,-K}_r(2,c_1,c_2)=\bigl\{E\in
M^H(r,c_1,c_2)\,\bigm|\,h^1(E)\ge-\chi(E)+r\bigr\}.
 \tag1.9
 $$
But for surfaces, the last inequality can be rewritten
 $$
h^1(E)\ge-\chi(E)+r\iff h^0(E)+h^2(E)\ge r.
 $$
Hence we have a decomposition
 $$
\sM^{g_H,-K}_r(2,c_1,c_2)=\bigcup_{i+j=r}M_{i,j}^H(2,c_1,c_2),
 $$
where the components on the right-hand side are the algebraic subvarieties of
$M^H(2,c_1,c_2)$ defined by
 $$
M_{i,j}^H(2,c_1,c_2)=\bigl\{E\in M^H(2,c_1,c_2)\,\bigm|\,h^0(E)\ge i,h^2(E)\ge
j\bigr\}.
 \tag1.10
 $$

On the other hand, sending $E\rightsquigarrow E^*(K)=E^*\otimes\Oh_S(K)$
identifies
 $$
M^H(2,c_1,c_2)=M^H(2,rK-c_1,c_2-c_1\cdot K-K^2),
 $$
and Serre duality gives
 $$
M_{i,j}^H(2,c_1,c_2)=M_{j,i}^H(2,rK-c_1,c_2-c_1\cdot K-K^2).
 \tag1.11
 $$
Now the Gieseker compactification $\overline{M^H(2,c_1,c_2)}$ (see \cite{G})
gives the compactifications $\overline{M_{i,j}^H(r,c_1,c_2)}$.

The standard definition of $\mu$@-homomorphism in the algebraic geometric
context
(see \cite{T1} or \cite{O'G}) gives the collection of polynomials
 $$
a \ga^H_{i,j}(2,c_1,c_2),
 \tag1.12
 $$
Now to compute the algebraic geometric version of the spin polynomial of
jumping level $r$ (1.6), we must sum the individual polynomials (1.12)
 $$
a\ga^H_r(2,c_1,c_2)=\sum_{i+j=r}a\ga^H_{i,j}(2,c_1,c_2).
 \tag1.13
 $$
But here care is needed, because the natural orientations of the components
(1.10) can be different (see \cite{P--T}, Chap\.~1, \S5).

One prove that the algebraic geometric polynomials (1.13) and spin polynomials
(1.6) are equal using the same arguments as for the original Donaldson
polynomials (see Morgan \cite{M}). More precisely, for a Hodge metric $g_H$, if
all the moduli spaces (1.10) have the expected dimension and avoid the
reducible
connections then
 $$
a\ga^H_r(2,c_1,c_2)=s_r\ga^{p_1}_{g_H,c_1-K_S};
 \tag1.14
 $$
see (1.13) and (1.7) and \cite{T3}, \S6.

To finish this important example we must recall the description of the coupled
Dirac operator (1.1) in terms of the complex structure of $S$. Namely
 $$
D^{g_H,-K_S,\nao}_{a=E}\:\Om^{0,0}(E)\oplus\Om^{0,2}(E)
 @>d'' \oplus(d'')^*>> \Om^{0,1}(E)
 \tag1.15
 $$
is the convoluted Dolbeault complex of $E$.

In the following section, which is the core of this paper we describe the
dependence of the spin polynomials on the choice of the Riemannian metric $g$.

\endexample

\head 2. The dependence on the metric \endhead

Let $\sS(M)$ be the space of Riemannian metrics on $M$. Every Riemannian metric
$g$ provides

1) an identification
 $$
H^2(M,\R)=H_g
 $$
where $H_g$ is the space of harmonic 2-forms with respect to the metric $g$.

2) a Hodge $*$ decomposition
 $$
H_g=H_g^+\oplus H_g^-
 $$
where $H_g^+$ is the subspace of selfdual forms and $H_g^-$ that of
antiselfdual
forms. Fixing an orientation of each $H_g^+$, we get the {\it period map}:
 $$
\Pi\:\sS(M)\to\Om\subset\Grass_{\text{orient}}(b_2^+,H^2(M,\R))
\qquad\text{defined by}\qquad\Pi(g)=H_g^+;
 \tag2.1
 $$
here $\Om$ is the subdomain of the Grassmanian of oriented subspaces in
$H^2(M,\R)$ of rank $b_2^+$ defined by:
 $$
\Om=\bigl\{H\in\Grass_{\text{orient}}(b_2^+,H^2(M,\R))\,\bigm|\,q_M{}\rest
H>0\bigr\},
 \tag2.2
 $$
where $q_M$ is the intersection form of $M$.

Any class $e\in H^2(M,\Z)$ with $e^2<0$ defines a subspace in $\Om$:
 $$
W_e=\bigl\{H\in\Om\,\bigm|\,e\in H^{\perp}\bigr\}.
 \tag2.3
 $$

A metric $g\in\sS$ is called $e$@-{\it irregular} if
 $$
\Pi(g)\in W_e,
 \tag2.4
 $$
and the corresponding subset of metrics
 $$
\sS_e=\bigl\{g\in\sS\,\bigm|\,\Pi(g)\in W_e\bigr\}\subset\sS
 \tag2.5
 $$
is called a {\it hurdle}; a hurdle $\sS_e$ is called a {\it wall\/} if its
complement
 $$
\sS\setminus\sS_e
 \tag2.6
 $$
is disconnected.

The topological type $(c_1\mod2,p_1)$ of a $\SO(3)$@-bundle $E$ on $M$ defines
a system $\{e\}_E$ of vectors given by the conditions:
 $$
\aligned
1)\qquad&e=c_1\mod2\\
2)\qquad&e^2=p_1\\
\endaligned
 \tag2.7
 $$
and a system of hurdles:
 $$
\{\sS_e\},\qquad\text{for}\quad e\in\{e\}_E.
 \tag2.8
 $$

Let $e\in\{e\}_E$ and suppose that $g$ is an $e$@-irregular metric. What does
this mean geometrically? First of all, by condition 1) of (2.7),
 $$
e=c_1-2\de,
 \tag2.7$'$
 $$
and $E$ splits topologically as a sum of line bundles:
 $$
E=L_{\de}\oplus L_{c_1-\de}
 \tag2.9
 $$

\remark{Remark} The choice of $\de $ or $c_1-\de$ defines an orientation
of the hurdle $\sS_{c_1-2\de}$ because of the equality
 $$
-e=c_1-2(c_1-\de).
 $$
So it is convenient to define an oriented hurdle by the class $\de$.
\endremark

Now if the line bundle $L_{\de}$ has a $\U(1)$@-connection $\al$ and
$L_{c_1-2\de}$ a $\U(1)$@-connection $\al'$ such that the induced connection
 $$
\al'\otimes\al^*\quad\text{on}\quad L_{c_1-\de}\otimes L_{\de}^*\subset\ad E
 $$
is $g$@-antiselfdual, then the reducible connection $\al\oplus\al'$ on $E$
defines a singular point
 $$
\al\oplus\al'\in\Sing\sM^g(E).
 \tag2.10
 $$
Thus a metric $g$ is $e$@-irregular if and only if the instanton space
$\sM^g(E)$ admits a singular point of type (2.9).

Now from (2.3) it is easy to see that $W_e$ has codimension
 $$
\codim W_e=b_2^+
 $$
in the domain $\Om$ (2.2). From this and the transversality conditions of the
period map (2.1) along $W_e$, one can prove that if $b_2^+>1$ then a hurdle is
{\it never} a wall, that is, the complement
 $$
\sS\setminus\bigcup_{e\in\{e\}_E}\sS_e
 $$
is connected.

One deduces from this, using the bordism arguments of \cite{D--K}, that if
$b_2^+>0$ then the Donaldson and spin polynomials of any jumping level are
independent of the metric, assumed to be regular in the sense of Uhlenbeck.

Now we would like to investigate when the singular point (2.10) is contained in
the subspace $\sM_r^{g,C}(E)$ of the instanton space $\sM_g(E)$ (see (1.4)).

By definition, if
 $$
\al\oplus\al'\in\sM_r^{g,C}(E)
 \tag2.11
 $$
then the operators
 $$
\aligned
 &D^{g,C}_{\al}\:\Ga^\infty(L_{\de}\otimes W^+)\to\Ga^\infty
(L_{\de}\otimes W^-)\\
 &D^{g,C}_{\al'}\:\Ga^\infty(L_{c_1-\de}\otimes W^+)\to\Ga^\infty
(L_{c_1-\de}\otimes W^-),\\
\endaligned
 \tag2.12
 $$
have kernels satisfying
 $$
\rank\ker D^{g,C}_{\al}+\rank\ker D^{g,C}_{\al'}\ge r.
 \tag2.13
 $$

Now a metric $g \in \sS$ is called {\it Dirac regular}, if
 $$
\rank\ker D^{g,C}_{\al}\ge r\implies\chi_C(L)\ge r
 \tag2.14
 $$
for the $\U(1)$@-connection $\al$ on any line bundle $L$; here $\chi_C(L)$
is the index of the couple Dirac operator (2.12). It is proved that a generic
metric $g\in\sS$ is Dirac regular (see for example \cite{P--T}, Chap\.~1, \S1).

It is convenient to define a {\it $\Spin^C$ hurdle} in terms of the class $\de$
(2.7$'$) and the decomposition (2.9):
 $$
\sS_{\de}=
\bigl\{g\in\sS_{c_1-2\de}\,\bigm|\,\text{$g$ is Dirac irregular}\bigr\}.
 \tag2.15
 $$
The main observation is the following
 $$
\text{a $\Spin^\C$ hurdle $\sS_{\de}$ is a wall
if and only if $\chi_C(L_{\de})\ge r$.}
 $$
Informally, this is easy, because both the $(c_1-2\de)$@-irregularity condition
and the Dirac irregularity condition are conditions of codimension $\ge1$. If
we
can prove that they are independent conditions then we are done.

To carry out these arguments rigorously, we have to use a very simple trick. We
consider a new space of parameters for our family of operators, larger than the
space of all Riemannian metrics $\sS$, namely the direct product
 $$
\sS\times\Om^2;
 \tag2.16
 $$
that is, the set of pairs $(g,\nao)$ consisting of any metric $g$ and any
Hermitian connection $\nao$ on the determinant spin bundle $L_C$ (see (1.1)),
which we can view as a 2-form (for a fixed metric $g$).

Of course, there is no advantage for the lifted walls of a
$(c_1-2\de)$@-irregular metric, but we can use the new parameter $\nao$ to
regularise the Dirac irregularity. Namely on twisting the spinor bundles
$W^\pm$ by the line bundle $L_{\de}^*$, and changing the $\Spin^\C$ structure
from $C$ to $C'=C-2\de$, the coupled Dirac operator $D^{g,C}_{\al}$ (2.12)
becomes the ordinary Dirac operator
 $$
D^{g,C-2\de}\:\Ga^\infty(W^+)\to\Ga^\infty(W^-)
 $$
of the metric $g$ (and a new $\Spin^\C$ structure). So we must study how the
jumping behaviour of the kernel of the ordinary Dirac operator of any
$\Spin^\C$
structure changes on deforming the metric $g$. But this was done in
\cite{P--T},
Formulas (1.3.3) and (1.3.4), Proposition 3.1.1 and Corollary.

{}From this, the description of the normal cone to $W_{c_1-2\de}$ and the
conditions of elliptic regularity provide the transversality conditions along a
$\Spin^\C$ hurdle. Now we can specify conditions for a wall for the moduli
space
of jumping instantons of type $(r,C+c_1,p_1)$ (see (1.7)): a class
$\de\in H^2(M,\Z)$ defines a wall if and only if
 $$
\aligned
1)\qquad&0\ge(c_1-2\de)^2\ge p_1\\
2)\qquad&\chi_C(L_{\de})\ge r\qquad
\text{(respectively, $\chi_C(L_{(c_1-\de)})\ge r$)}
\endaligned
 \tag2.17
 $$

\remark{Remark 1} Why do we only consider the condition on one component of the
decomposition (2.9)? The reason is the following: consider only the case
$\chi_C(E)\le0$. If one component $L_{\de}$ has $\chi_C(L_{\de})>0$ then the
other component $L_{c_1-\de}$ has $\chi_C(L_{c_1-\de})<0$.
\endremark

\remark{Remark 2} Thus the crucial difference between the space of instantons
$\sM^g$ and the subspace $\sM^{g,C}_r$ (1.4) is that $\sM^g$ depends only on
the
conformal class of $g$ whereas the subspace of jumping instantons has as its
full collection of indexes $\sM^{g,C,\nao}_r$ (as usual we drop the final
index).
\endremark

The conditions 1) and 2) of (2.17) imply the following two inequalities:
 $$
\aligned
1)\qquad&\frac14c_1^2\le c_1\cdot\de-\de^2\\
2)\qquad&\chi_C(L_{\de})=\frac12\de(\de+C)+\frac18(C^2-I)\ge r\\
\endaligned
 \tag2.18
 $$
where $I$ is the index of $M$. (This is the Atiyah--Singer formula.)

The second inequality is equivalent to the following:
 $$
-C\cdot\de-\de^2\le\frac14C^2-\frac14I-2r
 \tag2.18$'$
 $$

\example{Very important case}
 $$
c_1=-C
 \tag2.19
 $$
Then the right-hand side of the inequality 1) is equal to the left-hand side of
the inequality (2.18$'$) and we have the inequality
 $$
8r\le-I
 \tag2.20
 $$
Hence {\it if $8r>-I$ then the system of walls is empty}.
\endexample

\remark{Remark} As the reader can probably guess, this paper is written for the
sake of the preceding sentence. From it, using the bordism arguments of
\cite{D--K} one obtains that {\it the spin polynomials of jumping level $8r>-I$
are independent of the metric}, assumed to be regular in the sense of
Uhlenbeck.

So the spin polynomial (1.7) with $C+c_1=0$, that is, $s_r\ga^{p_1}_{g,0}$ is
called the {\it spin canonical polynomial of jumping level $r$}.

These polynomials behave naturally under diffeomorphisms of $M$. Namely, if
$8r>-I$, then for any $\si\in H^2(M,\Z)$ and any $\fie\in\Diff M$ that
preserves
the orientation of a maximal positive subspace of $H^2(M,\R)$, we have
 $$
s_r\ga^{p_1}_{g,0}(\si)=s_r\ga^{p_1}_{g,0}(\fie(\si)).
 \tag2.21
 $$
This means that some aspects of the shape of these polynomials (and their
coefficients, of course) are invariants of the smooth structure of 4-manifolds.
\endremark

\remark{Remark} Of course, the same formula holds for every class $(C+c_1)$
which
is invariant under the group of diffeomorphisms $\Diff M$ preserving the
orientation of a maximal positive subspace of $H^2(M,\R)$. But it can happen,
at
least a priori, that the only invariant class is the trivial class 0.
\endremark

To explain why the spin polynomials (1.7) are {\it canonical}, we return to our
main example.

\example{Example. Algebraic surfaces} If $M$ is the underlying manifold of an
algebraic surface $S$, then the anticanonical class $-K_S$ gives a canonical
$\Spin^\C$ structure, and the Very Important Case (2.19) predicts the equality
 $$
c_1=K_S.
 \tag2.22
 $$
Then the standard numerical invariants of $S$, the topological Euler
characteristic $c_2(S)$, and $K_S^2$ satisfy
 $$
\frac13(2c_2(S)-K_S^2)=-I.
 $$
Moreover, by the Noether formula,
 $$
-I+K_S^2=8\bigl(\frac1{12}(K_S^2+c_2(S))\bigr)=8(p_g+1),
 $$
where $p_g$ is the geometric genus of $S$.

Hence the inequality $8r>-I$ is equivalent to the inequality
 $$
-K_S^2<8(r-1)-8p_g.
 $$
Of course, we are interested in the case $p_g=0$. Then our inequality is
 $$
-K_S^2<8(r-1) .
 \tag2.23
 $$
Therefore, for a minimal surface of general type, jumping level 1 is already
enough to guarantee the invariance of the spin canonical polynomials.
\endexample

Returning to the general case, it is convenient to write the pair of conditions
(2.17) in the equivalent form:
 $$
\aligned
1)\qquad&0\ge(c_1-2\de)^2\ge p_1,\\
2)\qquad&(C+2\de)^2\ge 8r+I,\\
\endaligned
 \tag2.17$'$
 $$
and to consider a new vector
 $$
\De=-C-2\de.
 \tag2.24
 $$
Then the pair of conditions (2.17) is equivalent to the pair
 $$
\aligned
1)\qquad&0\ge(\De+(c_1+C))^2\ge p_1,\\
2)\qquad&\De^2\ge8r+I.\\
\endaligned
 \tag2.25
 $$

Now to convince you that spin canonical polynomials are useful in studying the
geometry of differentiable 4-manifolds we have to prove that these invariants
are
nonvanishing. As usual, we have to restrict ourselves to our main example, the
case of algebraic surfaces. Recall that Donaldson and Zuo proved that if
$|p_1|$
is large then any Hodge metric $g_H$ on algebraic surface is generic in the
sense of Freed and Uhlenbeck for the moduli space $\sM^g(c_1\mod2,p_1)$ (1.4).
In
spite of this, it is not true that every Hodge metric $g_H$ on an algebraic
surface is regular for the moduli space of jumping instantons. Fortunately we
can
describe what we need to add to the standard arguments of Donaldson theory to
use the new invariants. We do this in the following section.

\head 3. How much is the independence? \endhead

Return to the Very Important Case of algebraic geometry (see (1.8--14), and
(2.22--23)). This case is singular for the ad hoc reason that, by Serre
duality,
 $$
 M_{i,j}^H (2,K_S,c_2)=M_{j,i}^H (2,K_S,c_2)
 \tag3.1
 $$
(see (1.19--11)). This means that even if the space $M_{i,j}^H(2,K_S,c_2)$
has the right dimension, this locus has nontrivial multiple structure; that is,
as a subscheme, it is not reduced, it has nilpotents. Fortunately, we can
describe the scheme-theoretic structure precisely. We do this in the simplest
case $r=1$, that is, for the minimal jumping level.

First of all, recall that in the regular case, the fibre of the normal bundle
to
the sublocus $M_{1,0}^H(2,c_1,c_2)$ in $M^H(2,c_1,c_2)$ at a point $E\in
M_{1,0}^H(2,c_1,c_2)$ is given by
 $$
(N_{M_{1,0}\subset M})_E=\Hom(H^0(E),H^1(E)).
 \tag3.2
 $$
So the codimension of $M_{1,0}^H(2,c_1,c_2)$ in $M^H(2,c_1,c_2)$ is given by
 $$
\vcodim M_{1,0}^H(2,c_1,c_2)=h^0(E)(h^0(E)-\chi(E)).
 \tag3.3
 $$
Moreover if $E\in M_{1,0}^H(2,c_1,c_2)$ is a singular point then in the space
(3.2), we have the fibre of the {\it normal cone} of $M_{1,0}^H(2,c_1,c_2)$ in
$M^H(2,c_1,c_2)$:
 $$
(C_{M_{1,0}\subset M})_E\subset\Hom(H^0(E),H^1(E)),
 \tag3.4
 $$
defined as in \cite{F}.

Now in our Very Important Case, the spaces $H^i(E)$ involved in (3.2--4)
admit additional structures. Namely, by Serre duality
 $$
H^0(E)=H^2(E)^*\qquad\text{and}\qquad H^1(E)=H^1(E)^*.
 $$
This means that the vector space $H^1(E)$ has a nondegenerate {\it symmetric}
quadratic form
 $$
q_E\:H^1(E)\to H^1(E)^*,
 \tag3.5
 $$
and a {\it light cone} of isotropic vectors,
 $$
Q_E\subset H^1(E).
 \tag3.6
 $$

Consider the simplest case
 $$
\rank H^0(E)=h^0(E)=1.
 \tag3.7
 $$
A formal normal vector $n$ in the fibre of the formal bundle
$\Hom(H^0(E),H^1(E))$ is given as a nontrivial homomorphism
 $$
n\:H^0(E)\to H^1(E).
 \tag3.8
 $$
Then {\it $n$ is contained in the fibre of normal cone (3.4) if and only if the
image of the homomorphism (3.8) is contained in $Q_E$:}
 $$
n(H^0(E))\subset Q_E;
 \tag3.9
 $$
that is, the image is isotropic with respect to the quadratic form $q_E$. So
the
normal cone (3.4) is given by
 $$
(C_{M_{1,0}\subset M})_E=Q_E
 \tag3.9$'$
 $$
This is almost obvious: the Dirac operator for a Hodge metric $g_H$ is the
convoluted Dolbeault complex of $E$:
 $$
D^{g_H,-K_S\nao}_{a=E}\:\Om^{0,0}(E)\oplus\Om^{0,2}(E)
@>d''\oplus(d'')^*>>\Om^{0,1}(E)
 $$
(see (1.15)). Now the arguments used to prove the transversality theorem in
\cite{P--T}, Chap\.~1, \S3 make mathematical sense of the symbol
 $$
\frac{\partial\;}{\partial g}(D^{gH+\ep g,-K_S,\nao}_a)
 \tag3.10
 $$
as the line variation of a coupled Dirac operator with a jumping kernel with
a fixed identification
 $$
\aligned
H^0(E)&{}=\ker\frac{\partial\;}{\partial g}(D^{gH+\ep g,-K_S,\nao}_a)=\C,\\
H^1(E)&{}=\coker\frac{\partial\;}{\partial g}(D^{gH+\ep g,-K_S,\nao}_a).
\endaligned
 \tag3.11
 $$

So we have the diagram
 $$
\spreadmatrixlines{6pt}
\matrix
&\Om^{0,0}\kern-3cm&@>d''>>&\kern-3cm\Om^{0,1}&\\
&\uparrow\kern-3cm&&\kern-3cm\downarrow&\\
&H^0(E)=\ker\frac{\partial\;}{\partial g}(D^{gH+\ep g,-K_S,\nao}_a)&@>n>>&
\coker\frac{\partial\;}{\partial g}(D^{gH+\ep g,-K_S,\nao}_a)=H^1(E)&\\
\endmatrix
 \tag3.12
 $$
Now the Hermitian structure and Serre duality provide the diagram
 $$
\spreadmatrixlines{6pt}
\matrix
\Om^{0,0}(E)&@>d''>>&\Om^{0,1}(E)&@>d''>>&\Om^{0,2}\\
&&||&&||\\
&&(\Om^{0,1})^*&@<{(d'')^*}<<&(\Om^{0,2})^*
\endmatrix
 \tag3.13
 $$
But the Dolbeault complex is exact, hence the composite
 $$
H^0(E)@>n>>H^1(E)@>q_E>>H^1(E)^*@>n^*>>H^2(E)
 \tag3.14
 $$
is zero, where $q_E$ is the correlation (3.5). This means exactly (3.9).

\example{Example. The Barlow surface} Let $S$ be the Barlow surface (see
\cite{K} for the definitions and motivation) with
 $$
K_S^2=1,
 $$
and consider the moduli space
 $$
M^H(2,K_S,1).
 \tag3.15
 $$
Then this moduli space is a finite set of vector bundles, so is compact. More
precisely, the underlying reduced subscheme of $M^H(2,K_S,1)$ is given by
 $$
M^H(2,K_S,1)\red=\Bs|2K_S|,
 \tag3.16
 $$
the set of base points of the bicanonical pencil. More precisely if
 $$
\Bs|K_S|=(2K_S)^2=p_1+p_2+p_3+p_4,
 \tag3.17
 $$
is the quadruple of base points of $|K_S|$, then for every point $p_i$ there
exist, up to the action of $\C^*$, only one nontrivial extension of type
 $$
0\to\Oh_S\to E_i\to J_{P_i}(K_S)\to0,
 \tag3.18
 $$
and
 $$
M^H(2,K_S,1)\red=\{E_{i}\},\quad\text{for $i=1,\dots,4$.}
 \tag3.19
 $$
Now for every $i=1,\dots,4$, the Euler characteristic is $\chi(E_i)=1$,
and
 $$
h^0(E_i)=h^1(E_i)=h^2(E_i)=1.
 \tag3.20
 $$

As we know, $H^1(E_i)$ has a nontrivial quadratic form (see (3.5)), for which
zero is the only isotropic vector:
 $$
Q_{E_i}=0.
 $$

Consider the subscheme
 $$
M^H_{1,0}(2,K_S,1)\subset M^H(2,K_S,1).
 \tag3.21
 $$
Then
 $$
M^H_{1,0}(2,K_S,1)\red=M^H(2,K_S,1)\red=\{E_{i}\}.
 \tag3.22
 $$
Now using formula (3.9), we can see that
 $$
M^H_{1,0}(2,K_S,1)=M^H(2,K_S,1)
 \tag3.23
 $$
as schemes.

Indeed, a homomorphism $n\:H^0(E_i)\to H^1(E_i)$ is an element of the normal
cone of $M^H_{1,0}(2,K_S,1)$ in $M^H(2,K_S,1)$ if and only if $n=0$. But this
means exactly that the equality (3.23) holds. On the other hand, it means that
every $E_i$ has multiplicity 2 in $M^H_{1,0}(2,K_S,1)=M^H(2,K_S,1)$. So we get
a description of the nilpotent structure of $M^H_{1,0}(2,K_S,1)=M^H(2,K_S,1)$,
and we can see that under a generic deformation of the Hodge metric $g_H$ on
the underlying differentiable structure of $S$, the quadruple of points
(3.17)=(3.19) bifurcate to 8 regular instantons. So the Donaldson--Kotschick
polynomial of degree 0, which in this case equals half the spin canonical
polynomial, is given by
 $$
\ga^{E_i}_g=\ga^{K_S\mod2,-3}_{g_H}=8,\quad\text{and}\quad
s\ga^{-3}_{g_H,0}=16.
 \tag3.24
 $$
\endexample

Of course, we consider this example as an application of the description of the
scheme structure of the jumping instantons locus in the algebraic geometric
situation of the Very Important Case. In the following section we describe some
applications of spin canonical polynomials.

\head 4. Applications \endhead

We first recall the main constructions, results and conjectures concerning the
smooth classification of algebraic surfaces. Every compact nonsingular
algebraic
surface $S$ over $\C$ defines three underlying structures: its topological
class
$tS$, its underlying differentiable 4-manifold $dS$, and its deformation class
as an algebraic surface $vS$; compare the survey \cite{T1}, (5.44).

For any topological 4-manifold $X$, let $\diff(X)$ be the set of differentiable
4@-mani\-folds topologically equivalent to $X$. This set is discrete. Let
 $$
\alg_k(X)=\bigl\{M\in\diff(X)\,\bigm|\,M=dS,\ka(S)=k\bigr\}
 \tag4.1
 $$
be the subset of $\diff(X)$ containing the underlying structures of the
algebraic surfaces of Kodaira dimension $k$ and
 $$
\alg(X)=\bigcup_{k=-\infty}^2\alg_k(X)\subset\diff(X)
 \tag4.2
 $$
be the subset of the underlying differentiable structures of algebraic
surfaces.

\proclaim{Van de Ven conjecture}
 $$
k\ne k'\implies\alg_k\cap\alg_{k'}=\emptyset.
 \tag4.3
 $$
\endproclaim
This conjecture has recently been finally settled by Friedman and Qin
\cite{F--Q}, and independently by Pidstrigach cite{P}.

Finally, we set
 $$
\var(X)=\bigl\{vS\,\bigm|\,tS=X\bigr\}
 \tag4.4
 $$
for the set of the deformation classes of algebraic surfaces of topological
type $X$.

By definition there is a surjection
 $$
f\:\var(tS)\to\alg(tS),
 \tag4.5
 $$
and the main question is to describe the fibres of $f$,
 $$
f^{-1}(dS)\subset\var(tS).
 \tag4.6
 $$

We have some experience in estimating these sets:
 $$
\card\alg(t\plane)=1;
 \tag4.7
 $$
this is a Corollary of Yau's theorem, proved independently in \cite{T2}.

\remark{Remark} This result should not be confused with the smooth Poincar\'e
conjecture (see the Introduction), which says that
 $$
\card\diff(t\plane)=1.
 $$
\endremark

For the topological type $tK$ of the K3 surface $K$, one has
 $$
\card\alg_0(tK)=1,\qquad\text{but}\qquad\card\alg_1(tK)=\infty;
 $$
see Friedman and Morgan, \cite{F--M}. In the same vein, for every surface of
Kodaira dimension 1, we have
 $$
\card\alg_1(tS)=\infty.
 $$
Finally, for surfaces of general type,
 $$
\card\alg_2(tS)<\infty
 \tag4.8
 $$
(see for example \cite{7}).

The structures $(vS,dS,tS)$ corresponds to three groups
 $$
\Mon S\subset\Mod S\subset\O(q_S)
 \tag4.9
 $$
where $\O(q_S) $ is the orthogonal group of the lattice $H^2(tS,\Z)$ with the
intersection form $q_S$ (recall that this quadratic form defines the topology
type $tS$ uniquely), $\Mod S$ is the image of the standard representation of
the
diffeomorphisms group of $dS$ preserving the orientation to $\O(q_S)$ and $\Mon
S$ is the subgroup of $O (q_S)$ generated by all monodromy automorphisms of all
algebraic families of the surface $S$. The algebraic classification of surfaces
is closely related to the differentiable classification of the underlying
4-manifolds (the fibre (4.6) is a measure of this relation) and the subgroup
$\Mon S\subset\Mod S$ can be large enouph to describe this relation in some
partial cases.

Our first application of the spin canonical polynomials is to (4.3).

\example{Application 1. Van de Ven conjecture} Following the results of
Friedman and Morgan \cite{F--M}, the only case of (4.3) that remains to prove
is
 $$
\alg_2\cap\alg_{-\infty}=\emptyset.
 \tag4.10
 $$

It has been observed many times (see \cite{T3} or \cite{T4}) that all the spin
canonical polynomials of a rational surfaces vanish. The reason for this is the
following: if a polarisation $H$ on $S$ satisfies
 $$
K_S\cdot H\le0
 $$
then a torsion free sheaf $F$ has either a section or a cosection $F\to K_S$
that contradicts the stability of $F$. Hence all the spaces (1.10) with
$c_1=K_S$ are empty, and all the spin canonical polynomials vanish.

On the other hand, for a surface of general type $S$, the map
 $$
m\:S\to S\mi
 $$
to its minimal model $S\mi$, and its collection of exceptional rational
curves
 $$
\bigl\{l_i\subset S\,\bigm|\,l_i\iso\proj^1,l_i^2=-1\bigr\},
\quad\text{for $i=1,\dots,n$}
 \tag4.11
 $$
are uniquely determined. Let $K\mi=m^*(K_{S\mi})$ be the pullback of the
canonical class of the minimal model.

Consider first the case $K_S^2>0$, that is, either $S$ is minimal or the number
of the blown points is less then $K\mi^2$. Then it is easy to see that a
general nontrivial extension of type
 $$
0\to\Oh_S\to E\to J_{\xi}(K\mi)\to0,
 \tag4.12
 $$
with $\xi\in\Hilb^dS$ a general zero-dimensional subscheme of large degree $d$,
is $H$@-stable for any polarization $H$, and
 $$
E\in M^H_{1,0}(2,K_S,d)\implies M^H_{1,0}(2,K_S,d)\ne\emptyset
 \tag4.13
 $$
We can prove this using the same arguments as in the proof of \cite{T3},
Theorem~4.1. Thus, if
 $$
\dim M^H_{1,0}(2,K_S,d)=\vdim M^H_{1,0}(2,K_S,d)
 \tag4.14
 $$
then the spin canonical polynomial is nonzero:
 $$
s\ga^{K_S^2-4d}_{g_H,0}\ne0,
 \tag4.15
 $$
and spin canonical polynomials of jumping level 1 distinguish the smooth types
of
rational surfaces from surfaces of general type.

If equality doesn't hold in (4.14) then we can use the regularisation procedure
described in detail in \cite{T3}, \S7, and obtain the inequality (4.15), and
this
again distinguishes a surface of general type $S$ from rational surfaces.

On the other hand, if $K_S^2<0$ then for simple arithmetical reasons there
exists a vector $\de$ with
 $$
\de\cdot(\de-K_S)\ge2(r-1)
\quad\text{and}\quad
-K_S^2<8(r-1),
 \tag4.16
 $$
and a polarization of the form
 $$
H=H\mi-\sum_{i=1}^{n}a_il_i,\quad\text{with $a_i>0$,}
 \tag4.17
 $$
where $H\mi$ is the inverse image of an almost canonical polarization on
$S\mi$, such that
 $$
2\de\cdot H\mi<K_S\cdot H\mi.
 \tag4.18
 $$
Recall (see \cite{T5}) that a polarisation $H_0$ is {\it almost canonical} if
the ray $\R^+\cdot H$ in the projectivised K\"ahler cone of $S\mi$ is close
to the ray $\R^+\cdot K\mi$ in the Lobachevsky metric.

Then it is again easy to see that, for some suitable choice of degree $d$, and
$\xi\in\Hilb^dS$ a general zero-dimensional subscheme of degree $d$,
a general nontrivial extension of type
 $$
0\to\Oh_S(\de)\to E\to\sJ_\xi(K_S-\de)\to0,
 \tag4.19
 $$
is $H$@-stable for the polarization (4.17), and
 $$
E\in M^H_{r,0}(2,K_S,d)\implies M^H_{r,0}(2,K_S,d)\ne\emptyset,
 $$
in the same vein as (4.12). So if
 $$
\dim M^H_{r,0}(2,K_S,d)=\vdim M^H_{r,0}(2,K_S,d),
 \tag4.14$'$
 $$
then the spin canonical polynomial is nonzero:
 $$
s_r\ga^{K_S^2-4d}_{g_H,0}\ne0,
 \tag4.15$'$
 $$
and the spin canonical polynomial of jumping level $r>\frac18(-K_S^2)+1$
distinguishes a surface of general type $S$ from rational surfaces.

If the equality (4.14$'$) doesn't hold, then we could use the regularisation
again, although I have not carried this out in detail.
\endexample

\example{Application 2. The reducibility of $\Mod S$}

We will consider $\Mod S$ as a subgroup of the orthogonal group $O(q_S)$, and
its representation as a transformation group of $H^2(S,\Z)$. We proved in
\cite{T4} that this representation is reducible: if $p_g>0$, there exist a
proper sublattice $sV(S)$ invariant under diffeomorphisms. More precisely,
 $$
\Z\cdot K_S\subset sV(S)\subset\Span{K_S,C_1,\dots,C_N},
 \tag4.20
 $$
that is $sV(S)$ contains $K_S$, and is contained in the lattice generated by
$K_S$ and the classes of all effective curves $\{C_i\}$ satisfying the
inequalities
 $$
2C_i\cdot K\mi\le K\mi^2
 \tag4.21
 $$
These classes are algebraic; hence if $p_g>0$ then $sV(S)$ is a proper
sublattice.

Now the main question is:
 $$
\text{how close is the sublattice $\Z\cdot K_S$ to $sV(S)$?}
 $$

Other diffeomorphism invariant sublattices are known: perhaps the best
approximation to $\Z\cdot K_S$ is the Kronheimer--Mrowka--Witten sublattice
 $$
L_{\text{KMW}}=\Span{C_1,\dots,C_k}
 \tag4.21
 $$
generated by the irreducible components of a general curve of the canonical
linear system $|K_S|$ (see \cite{W}). It is easy to see that
$L_{\text{KMW}}\subset sV(S)$

The method of proving diffeomorphism invariance is the following: the first
step
is to see the {\it shape} of the invariant polynomial. For example, if $p_g>0$
then Kronheimer, Mrowka and essentially Witten proved that the Donaldson
polynomials
 $$
\ga^{w_2(S),p_1\ll0}\in S^d(q_S,C_1,\dots,C_k)
 \tag4.22
 $$
are contained in the subring of $S^*H^2(S,\Z)$ generated by the intersection
form
as a quadratic form and the classes $C_i$ as linear forms. From this, the
invariance of the polynomials implies the invariance of the sublattice
generated
by the collection of linear forms. To see that the same thing holds for the
spin
sublattice $sV(S)$ (4.20) is much easier. So here we would like to explain why
the spin canonical polynomials belong to the subring
 $$
s\ga^{p_1\ll0}_0 \in S^d (q_S, C_1,\dots,C_N),
 \tag4.23
 $$
where $\{C_i\}$ is the collection of classes (4.21). To do this informally, it
is
very convenient to use the {\it vortex equation} and {\it the moduli space of
stable pairs} (see the introduction to \cite{B--D}).

Recall that for any $\U(2)$@-vector bundle $E$ on $S$ (or on any compact
K\"ahler
mani\-fold) there is a {\it Yang--Mills--Higgs functional\/} on the space of
pairs $(a,\fie)$, where $a$ is a $\U(2)$@-connection on $E$ and $\fie$ is any
section of $E$, depending on a real parameter $\tau$:
 $$
\text{YMH}_{\tau}(a,\fie)=\Vert
F_a\Vert^2+\Vert\dd_a''\fie\Vert^2+\frac14\Vert|\fie|_h^2-\tau\Vert^2
 $$
The pair $(a,\fie)$ is an absolute minimum of this functional if
 $$
\dd_a''\fie=0\quad\text{and}\quad
\bigwedge F_a-\frac i2|\fie|^2_h+\frac i2\tau=0
 \tag4.24
 $$
This system of differential equations is called the $\tau$@-{\it vortex
equation}, and the space $\sM_{\tau}$ of solutions up to gauge equivalence is
called the {\it moduli spaces of $\tau$@-vortices}.

For every class $\si\in H^2(S,\Z)$, let
 $$
\deg\si=\si\cdot[\om]
 $$
where $\om$ is the K\"ahler form.

Now it is convenient to change the parameter $\tau$:
 $$
\si=\frac1{4\pi}\tau\Vol S-\frac1{2}\deg E,
 \tag4.25
 $$
(where $\deg E=\deg c_1(E)$). Then Bradlow's Identification Theorem says that
 $$
\sM_{\tau}=\text{MP}_{\si}
 \tag4.26
 $$
is the {\it moduli space of $\si$@-stable pairs} $(E,s)$, where $E$ is a
holomorphic bundle and $s$ is a holomorphic section of $E$. Recall that a pair
$(E,s)$ is {\it $\si$@-stable} if
 $$
\deg L<
\cases
\half\deg E-\si &\text{if $s\in H^0(L)$,}\\
\half\deg E+\si &\text{otherwise,}\\
\endcases
\qquad\text{for any line subbundle $L\subset E$.}
 \tag4.27
 $$

Each moduli space $\text{MP}_{\si}$ is a family of vector bundles, and hence
slant product defines a cohomological correspondence
 $$
\mu\:H_i(S,\Z)\to H^{4-i}(\text{MP}_{\si}),
 $$
and a collection of polynomials
 $$
\ga P_{\si}\in S^*H^2(S,\Z).
 \tag4.28
 $$
\endexample

\remark{Remark} Actually, to define these polynomials we must either construct
the compactification of the moduli spaces $\text{MP}_{\si}$ or use a trick due
to Donaldson as in \cite{T3}, Lemma~2.1.
\endremark

Now from the definition it is easy to see that
 $$
\text{MP}_{\si}\ne\emptyset\implies\half\deg E\ge\si\ge0,
 \tag4.29
 $$
and for obvious numerical reasons the $\si$@-stability condition (4.27) remains
the same (and implies proper stability) for any
 $$
\si\in(\max(0,\half\deg E-i-1),\half\deg E-i).
 \tag4.30
 $$
So for $\si$ in this interval we get a fixed moduli space $\text{MP}_{\si}$,
and
as $\si$ varies we have a chain of birational transformations or flips (see
\cite{R})
 $$
\text{MP}_{\max}\leftrightarrow\dots\leftrightarrow\text{MP}_1
\leftrightarrow\text{MP}_0.
 \tag4.31
 $$
You can recognise the situation described by Thaddeus and Bertram for the
$\si$@-scale of the moduli space of stable pairs on an algebraic curve $C$ of
genus $g$ (see \cite{Th}, \cite{B--D} and \cite{B}). Recall that, in this {\it
classical} case, the moduli space on the extreme left has a birational regular
map
 $$
\text{MP}_{\max}\to M_{C,\xi}
 $$
to the moduli space of stable bundle with fixed determinant $\xi$ of degree
$2g-1$. In our case we have a map of the same type (now $c_1(E)=K_S$)
 $$
\text{MP}_{\max}\to M^H_{1,0}(2,K_S,d),
 \tag4.32
 $$
and from this it is easy to see that the polynomial (4.28)
 $$
\ga P_{\max}=s\ga_{0}^{K_S^2-4d}
 \tag4.33
 $$
is our spin canonical polynomial of jumping level 1.

Now in the classical case, the moduli space $\text{MP}_0$ on the right-hand
end of (4.31) is a projective space of dimension $3g-3$.
 $$
\text{MP}_0=\proj^{3g-3}.
 $$

In our case $\text{MP}_0$ is not so simple. Namely,
 $$
\text{MP}_0=\text{GAM}
 $$
is the space of all nontrivial extensions of type
 $$
0\to\Oh_S\to E\to\sJ_\xi(K_S)\to0
 \tag4.34
 $$
for all $\xi\in\Hilb^d S$, up to the action of $\C^*$ (see \cite{T3, T4, T4}).
For large $d$, this space is birational to the direct product
 $$
\Hilb^d S\times\proj^n.
 $$
It is not hard to compute the polynomial $\ga P_{0}$, and to see that
 $$
\ga P_{0}\in S^{\cdot}(q_S,K_S).
 \tag4.35
 $$
Now we can compute the increment
 $$
\ga P_{1}-\ga P_{0},
 \tag4.36
 $$
For this, we have to blow up the space $\text{MP}_0$ (4.34) in the subvariety
of
all extensions admitting a diagram of type
 $$
\spreadmatrixlines{6pt}
\matrix
&&0\\
&&\uparrow\\
&&\sJ_\eta(C_1)\\
&^{\fie'}\nearrow&\uparrow\\
0\to\Oh&\to&E&\to&\sJ_\xi(K_S)\to0,\\
&&\uparrow&{}^\fie\nearrow& \\
&&\kern-1cm\Oh(K_S-C_1)\kern-1cm\\
&&\uparrow\\
&&0\\
\endmatrix
 \tag4.37
 $$
where $C_1$ is a curve of degree 1, and to carry out the {\it elementary
transformation} of the universal extension as in \cite{T4}, \S6. Computing the
slant product again for the new family of vector bundles, we get the shape of
the
increment (4.36)
 $$
\ga P_{1}-\ga P_{0}\in S^*(q_S,K_S,C_1).
 $$
Therefore
 $$
\ga P_{1}\in S^*(q_S,K_S,\{C_1\}),
 \tag4.38
 $$
where $\{C_1\}$ is the set of all curves of degree 1.

Continuing this procedure gives the shape of the spin canonical polynomial
(4.33).

In fact, we can't do this construction rigorously because there isn't any
rigorous treatment of the compactification of moduli spaces of $\si$@-stable
pairs on a K\"ahler surface, and there are other technical problems. But in the
classical case, there is the pure algebraic geometric procedure imitating this
K\"ahler geometric procedure proposed by Bertram in \cite{Be}. We can do the
same for an algebraic surface (see \cite{T4}).

Moreover, Huybrechts and Lehn in \cite{H--L} proposed a beautiful algebraic
geometric theory of stable pairs on a surface. First of all, they construct
compactifications of moduli spaces of stable pairs by torsion free sheaves.
After this, in place of the real number $\si$ (4.25), it is very natural to
consider a polynomial $\de(z)$ with rational coefficients such that $\de>0$ for
all $z\gg0$ (in fact $\si$ is the first coefficient of $\de$). Now we have an
exact order in the sequence of flips (4.31), and the elementary transformations
of the universal sheaves such as (3.37).

Our last remark is the following: we can prove the reducibility of $\Mod S$
even
if $p_g=0$, provided that $K_S^2>0$ and the sublattice $sV(S)$ (4.20) is
proper.
For this, it is enouph to verify that the sublattice $\Span{K_S,C_1,\dots,C_N}$
(see (4.20--21)) is proper.

\example{Example. The Barlow surface} The Barlow surface $S$ is minimal, and
has
four smooth rational $-2$@-curves $C_1,C_2,C_3,C_4$. Thus
 $$
C\cdot K_S\le K_S^2=1\implies C=\sum_{i=1}^4n_iC_i,\quad\text{with}\quad
n_i\ge0.
 $$
But $\rank\Pic S=9$, and therefore
 $$
\Span{K_S,C_1,C_2,C_3,C_4}\subset\Pic S=H^2(S,\Z)
 $$
is a proper sublattice.
\endexample

 \remark{Remark} The diffeomorphism invariant sublattice $sV(S)$ (4.20) is a
sublattice of $\Span{K_S,C_1,C_2,C_3,C_4}$, a priori, a proper sublattice. This
holds for the Barlow surface. Actually, the exact computation shows that
 $$
sV(S)=\Z\cdot K_S
 $$
(compare the computations for $-2$@-curves in \cite{P--T}, Chap.~4, \S3). Thus
the canonical class of the Barlow surface is diffeomorphism invariant. This
means that the following conjecture is now verified for all currently known
algebraic surfaces.
 \endremark

\proclaim{Conjecture} 1) Every simply connected surface of general type $S$ has
a proper sublattice $sV(S)\subset H^2(S,\Z)$ invariant under diffeomorphisms;
that is,
 $$
 dS'\in\diff(dS)\implies sV(S)=sV(S').
 $$
Moreover, $sV(S)$ satisfies
 $$
\align
&2)\qquad sV(S)\subset\bigcap_{dS'\in\diff(dS)}\Pic S';\\
\text{and}\quad&3)\qquad\bigcup_{dS'\in\diff (dS)}K_{S'}\subset sV(S).
\endalign
$$

\endproclaim


\Refs
 \widestnumber \key{EMG}

 \ref \key B
 \by W. Barth
 \paper Moduli of vector bundles on the projective plane
 \jour Invent. math.
 \vol 42
 \yr 1977
 \pages 63--91
 \endref

 \ref \key Be
 \by A. Bertram
 \paper Moduli of rank 2 vector bundles,theta divisors and the geometry of
curves in projective space.
 \jour J. Diff.Geom.
 \vol 35
 \yr 1992
 \pages 429--469 \endref

 \ref \key B--D
 \by S. Bradlow and G. Doskalopoulos
 \paper Moduli of stable pair for holomorphic bundles over riemann surfaces
 \jour Int.J.Math.
 \vol 2
 \yr 1991
 \pages 477--513 \endref

 \ref \key C
 \by F. Catanese
 \paper Connected components of moduli spaces
 \jour J. Diff.Geom.
 \vol 24
 \yr 1986
 \pages 395--399 \endref

 \ref \key D
 \by S. Donaldson
 \paper Polynomial invariants for smooth 4-manifolds
 \jour Topology
 \vol 29
 \yr 1990
 \pages 257--315 \endref

 \ref \key D--K
 \by S. Donaldson and P. Kronheimer
 \book The Geometry of Four-Manifolds
 \publ Clarendon Press
 \publaddr Oxford
 \yr 1990 \endref

 \ref \key F--U
 \by D. Freed and K. Uhlenbeck
 \book Instantons and four-manifolds
 \publ M.S.R.I.Publ. Springer
 \publaddr New York
 \yr 1988 \endref

 \ref \key F--M
 \by R. Friedman, J. Morgan
 \paper On the diffeomorphism type of certain algebraic surfaces
 \jour J. Diff. Geom.
 \vol 27
 \yr 1988
 \pages 371--398 \endref

 \ref \key F--Q
 \by R. Friedman, Z. Qin
 \paper On complex surfaces diffeomorphic to rational surfaces
 \jour File server alg-geom 9404010
 \yr 1994
 \pages (34 pp.) \endref

 \ref \key F
 \by W. Fulton
 \inbook Intersection Theory
 \publ Springer-Verlag
 \publaddr Timbuktoo
 \yr 1987 \endref

 \ref \key G
 \by D. Gieseker
 \paper On the moduli of vector bundles on an algebraic surface.
 \jour Ann. of Math.
 \vol 106
 \yr 1977
 \pages 45--60 \endref

 \ref \key O'G
 \by K.G. O'Grady
 \paper Algebro-geometric analogues of Donaldson polynomials
 \jour Inv. Math.
 \vol 107
 \yr 1992
 \pages 351--395 \endref

 \ref \key H
 \by N. Hitchin
 \paper Poncelet polygons and the Painlev\'e equations
 \jour Warwick Preprint 6
 \yr 1994
 \endref

 \ref\key H--L
 \by D. Huybrachts and M. Lehn
 \paper Stable pairs on curves and surfaces
 \jour Max Planck Inst\. Preprint MPI/92-93
 \yr 1993
 \endref

 \ref \key K
 \by D.Kotshick
 \paper On manifolds homeomorphic to $\plane\#8\overline\plane$
 \jour Invent. Math.
 \vol 95
 \yr 1989
 \pages 591--600 \endref

 \ref \key M
 \by J. Morgan
 \paper Comparison of the Donaldson polynomial invariants with their
algebro geometric analogs.
 \jour Princeton preprint
 \yr 1992 \endref

 \ref \key P
 \by V. Pidstrigach
 \paper Some glueing formulas for spin polynomials and the Van de Ven
conjecture
 \jour Izv. RAN (to appear)
 \vol 54
 \yr 1994
 \pages (28 pp.)
 \lang Russian
 \endref

 \ref \key P--T
 \by V. Pidstrigach and A. Tyurin
 \paper Invariants of the smooth structures of an algebraic surfaces arising
from the Dirac operator
 \jour Izv. AN SSSR
 \vol 52:2
 \yr 1992
 \pages 279--371
 \lang Russian
 \transl \nofrills English transl. in
 \jour Warwick preprint
 \yr 1992
 \vol 22 \endref

 \ref \key R
 \by M. Reid
 \paper What is a flip?
 \jour Preprint Utah
 \yr 1993 \endref

 \ref \key Th
 \by M. Thaddeus
 \paper Stable pairs,linear system and the Verlinde formula
 \jour Preprint
 \yr 1992 \endref

 \ref \key T1
 \by A. Tyurin
 \paper Algebraic geometric aspects of smooth structures. I, Donaldson
polynomials
 \jour Russian Math.Surveys
 \vol 44:3
 \yr 1989
 \pages 117--178 \endref

 \ref \key T2
 \by A. Tyurin
 \paper A simple method of distinguishing the underlying differentiable
structures of algebraic surfaces.
 \jour Mathematica Goettingensis,Sonderforschungsbereichs Geometry and Analysis
 \vol Heft 25
 \yr 1992
 \pages 1--24 \endref

 \ref \key T3
 \by A. Tyurin
 \paper Spin polynomial invariants of the smooth structures of algebraic
surfaces.
 \jour Mathematica Goettingensis, Sonderforschungsbereichs Geometry and
Analysis
 \vol Heft 6
 \yr 1993
 \pages 1--48
 \endref

 \ref \key T4
 \by A. Tyurin
 \paper Canonical and almost canonical spin polynomials of algebraic surfaces.
 \jour Proceedings of the Conference ''Vector bundles'', Durham. 1993.
 \vol Aspects of Mathematics.
 \yr 1993
 \pages 250--275 \endref

 \ref \key T5
 \by A. Tyurin \paper On almost canonical polynomials of algebraic surfaces.
 \jour Proceedings of the Conference ''Algebraic Geometry'', Yaroslavl. 1992.
 \yr 1993
 \pages 1--25 \endref

 \ref \key W
 \by E. Witten
 \paper Supersymmetric Yang--Mills theory on a four-manifolds.
 \jour Preprint IASSNS-HEP-94
 \yr 1994
 \endref

 \ref \key Z
 \by K. Zuo
 \paper Generic smoothness of the moduli of rank two stable bundles over an
algebraic surface
 \jour Preprint MPI
 \yr 1990
 \endref

 \endRefs
 \enddocument